 \documentclass[12pt]{amsart}
\usepackage{amsmath}
\usepackage{amsxtra}
\usepackage{amscd}
\usepackage{amsthm}
\usepackage{amsfonts}
\usepackage{amssymb}
\input amssym.def
\input amssym

\newtheorem{theorem}{Theorem}[section]
\newtheorem{corollary}{Corollary}[section]
\newtheorem{lemma}{Lemma}[section]

\newtheorem{remark}{Remark}[section]
\newtheorem{example}{Example}[section]

\newtheorem{proposition}{Proposition}[section]
\newcommand{\bea}{\begin{eqnarray}}
\newcommand{\eea}{\end{eqnarray}}
\newcommand{\Zp} {\Z _ {\ge 0} }

\newcommand{\ver}{L(-\frac{4}{3}\Lambda_0)}

\newcommand{\h}{{\frak h}}
\newcommand{\n}{{\frak n}}
\newcommand{\hn}{\hat {\frak n}}
\newcommand{\g}{{\frak g} }
\newcommand{\hg}{\hat { {\frak g } } }
\newcommand{\hh}{\hat {  {\frak h} } }

\def \la{\langle}
\def \a{\gamma  }

\newcommand{\ap}{{\alpha} }
\newcommand{\bp}{{\beta} }
\def \b{\delta  }
\def \ra{\rangle}
\newcommand{\C}{ \mathbb{C} }

\newcommand{\Z}{\mathbb{Z} }

\newcommand{\1}{\bf 1}
\def \l{\lambda}
\newcommand{\be}{\begin{equation}}
\newcommand{\ee}{\end{equation}}

\begin{document}

\title[]{A construction of admissible $A_1^{(1)}$--modules of level $-\frac{4}{3}$.}

\author{ Dra\v zen Adamovi\' c }

\date{}
\curraddr{Department of Mathematics, University of Zagreb,
Bijeni\v cka 30, 10 000 Zagreb, Croatia} \email {adamovic@math.hr}
\markboth{Dra\v zen Adamovi\' c} {}
\bibliographystyle{amsalpha}
  \maketitle

\begin{abstract}
By using    generalized vertex algebras associated to  rational
lattices, we construct explicitly the   admissible modules
for the affine Lie algebra $A_1 ^{(1)}$ of level $-\frac{4}{3}$. As an
application,  we show that the $\mathcal{W}(2,5)$ algebra with central
charge $c=-7$ investigated in \cite{A} is a subalgebra of the
simple vertex operator algebra ${\ver}$.
\end{abstract}

\section{Introduction}

In this paper we shall present an explicit construction of
irreducible representations of the vertex operator algebra
$L(-\frac{4}{3}\Lambda_0)$ associated to the irreducible vacuum
representation for the affine Lie algebra $A_1^{(1)}$ of level
$-\frac{4}{3}$.
 Explicit realizations of certain irreducible
highest weight representations play important role in the
representation theory of affine Lie algebras and the associated
vertex operator algebras.
  Integrable highest weight modules can be realized by
using the theory of lattice vertex algebras or Clifford vertex
superalgebras (cf. \cite{DL},\cite{F}, \cite{Fe},\cite{FFR},
\cite{K2}).
 Admissible modules are a broad class of irreducible highest
weight modules which contains integrable modules as a subclass.
Admissible modules have rational levels and their characters are
modular functions (cf. \cite{KW}, \cite{Wak}).
It was proved in \cite{AM} and \cite{DLM} that admissible $A_1
^{(1)}$--modules of level $k$ can be viewed as  modules for the
simple vertex operator algebra $L(k\Lambda_0)$. Moreover,  every
$L(k\Lambda_0)$--module from the category $\mathcal{O}$ is a
direct sum of admissible modules.

An explicit
construction of admissible  modules for the affine Lie algebra $C_n ^{(1)}$ of level $-\frac{1}{2}$
 was given in \cite{FF} by using Weyl
algebras. This construction was extended in \cite{We} to the
framework of vertex operator algebras. In the case $n=1$ this
construction gives the realization of four admissible
$A_1^{(1)}$--modules. It is important to notice that the vertex
algebras associated to  Weyl algebras can be realized by using
lattice vertex algebras (cf. \cite{FMS}, \cite{efren}, \cite{A2}).

In the physics literature the admissible representations are also
well-known. In particular, the Weyl vertex algebras correspond to
the $\beta \gamma$--system. Quite recently, the work \cite{LMRS}
presented a very detailed study of the $\beta \gamma$--systems and
its relation to the WZW--model $\widehat{su(2)}_{-\frac{1}{2}}$.
Their construction was based on   a $c=-2$  theory and a Lorenzian
boson.  M. Gaberdiel in \cite{Gab} presented an interesting  study
of the WZW--model $\widehat{su(2)}_{-\frac{4}{3}}$ and analyzed
the fusion products of admissible representations of level
$-\frac{4}{3}$. He also showed that this theory contains
logarithmic representations, i.e., the representations on which
the action of $L(0)$ is not diagonalisable. His approach was based
on the analysis of singular vectors and the concept of fusion
product.

 In the present paper we shall give an explicit
construction of admissible representations of level
$-\frac{4}{3}$. First we shall construct the simple vertex
operator algebra $L(-\frac{4}{3} \Lambda_0)$. Since the rank of
the vertex operator algebra ${\ver}$ is $-6$, it is natural to
consider vertex  algebras of rank $-7$ and rank one lattice vertex
algebras. Since $c=-7$ is the central charge of the $(1,p)$ models
for the Virasoro algebra (with $p=3$),  we can use  the
representation theory of the $(1,p)$  models. In \cite{A}, we
showed that the concept of generalized vertex algebras associated
to  rational lattices is useful for studying vertex algebras
associated to the $(1,p)$   models. It turns out that this
approach is useful for the present construction.

 Let us now describe
our construction in more details. We define the following lattices
\bea &&  \widetilde{L}= \Z \a + \Z \b, \quad \la \a , \a \ra = -
\la \b , \b \ra = \frac{1}{6}, \quad \la \a , \b \ra =0, \nonumber
\\ && L= {\Z} (\a + \b) + {\Z} (\a - \b). \nonumber \eea

Let $V_{\widetilde{L}}$ and $V_L$ be the associated  generalized
vertex algebras. We consider the subalgebra $V_{\Z \a}$ as a
module for the Virasoro algebra with central charge $-7$. As in
\cite{A}, we define the screening operator  $Q=e^{-6 \a} _0$.
 We prove, as one of
our main results, that the subalgebra of the generalized vertex
algebra $V_L$ generated by vectors
\bea
&& e = e^{3 (\a - \b)}, \quad  h  = 4 \b (-1) , \quad f  =
-\frac{2}{9} Q e^{3 (\a + \b)}, \nonumber
\eea
is isomorphic to the simple vertex operator algebra ${\ver}$ (for
details  see Section \ref{lattice-construction}). Then we consider
$V_L$ as a weak module for ${\ver}$. Although $V_L$ is not
completely reducible ${\ver}$--module, it contains some
irreducible submodules.  Since every ${\ver}$--module in the
category $\mathcal{O}$ is completely reducible (cf. \cite{AM}),
every highest weight vector in $V_L$ generates an irreducible
highest weight $A_1^{(1)}$--module. By using  this fact, we
identify admissible modules inside $V_L$. We also identify
infinitely many irreducible submodules of $V_L$ which are not
${\Zp}$--graded (cf. Theorem \ref{spectral}).

The  lattice construction enables us to investigate certain
subalgebras of the vertex operator algebra ${\ver}$. By using
structural results from \cite{A}, we will show that the simple
subalgebra $\mbox{Ker}_{{\ver}} h(0)$ of ${\ver}$ is isomorphic to
the tensor product of the $\mathcal{W}(2,5)$ algebra with central charge $-7$
and the free boson vertex algebra $M_{\b}(1)$ (cf. Theorem
\ref{izom-w}). As a consequence, we get that the coset vertex algebra
$C({\ver}, M_{\b} (1) )$ (= the coset $\frac
{\widehat{su(2)}_{-\frac{4}{3}} } {\widehat{u(1)} }$ ) is isomorphic to the $\mathcal{W}(2,5)$ algebra with central
charge $-7$.

\section{ Vertex operator algebra $L (k\Lambda_0)$ }

In this section we recall some basic facts about vertex operator algebras associated to affine Lie algebras
(cf. \cite{FZ}, \cite{Li-local}, \cite{MP}).

Let ${\g}$ be a finite-dimensional simple Lie algebra over ${\Bbb
C}$ and let $(\cdot,\cdot)$ be a nondegenerate symmetric bilinear
form on ${\g}$. Let  ${\g} = {\n}_- + {\h} + {\n}_+$    be a
triangular decomposition for ${\g}$.
 The affine Lie algebra ${\hg}$ associated
with ${\g}$ is defined as $ {\g} \otimes {\C}[t,t^{-1}] \oplus
{\C}c \oplus {\C}d, $ where $c$ is the canonical central element
\cite{K}
 and  the Lie algebra structure
is given by $$ [ x \otimes t^n, y \otimes t^m] = [x,y] \otimes t
^{n+m} + n (x,y) \delta_{n+m,0} c,$$ $$[d, x \otimes t^n] = n x
\otimes t^n $$ for $x,y \in {\g}$.  We   write $x(n)$ for $x
\otimes t^{n}$ and identify $\g$ with ${\g}\otimes t^{0}$.
The Cartan subalgebra ${\hh}$ and  the subalgebras ${\hg}_{\pm}$,
${\hn}_{\pm}$ of ${\hg}$ are defined as $${\hh} = {\h} \oplus
{\C}c \oplus {\C}d, \quad
 {\hg}_{\pm} =  {\g}\otimes t^{\pm1} {\C}[t^{\pm1}], \quad {\hn}_{\pm} =
 {\n}_{\pm} \oplus {\hg}_{\pm}.$$

Let  $P = {\g}\otimes {\C}[t] \oplus {\C}c \oplus {\C}d$ be upper
parabolic subalgebra.   For every $k \in {\C}$, $k \ne -h
^{\vee}$, let ${\C} v_k$  be $1$--dimensional $P$--module  such
that the subalgebra ${\g} \otimes  {\C}[t] + \C d$ acts trivially,
 and  the central element
$c$ acts as multiplication with $k \in {\C}$. Define the
generalized Verma module $N( k \Lambda_0)$ as
$$N({k}\Lambda_0) = U(\hg) \otimes _{ U(P) } {\C} v_k .$$
Then $ N(k \Lambda_0)$ has a natural structure of a vertex
operator algebra (VOA). The vacuum vector is ${\1} = 1 \otimes
v_k$.

Let $N^{1}(k \Lambda_0)$ be the maximal ideal in the vertex
operator algebra $N(k \Lambda_0)$. Then $L(k \Lambda_0) =
\frac{N(k \Lambda_0)}{N^{1}( k \Lambda_0)}$ is a simple vertex
operator algebra.

Let $M$ be any $U({\hg})$--module. A weight vector $v \in M$ is
called a singular vector if ${\hn}_{+} v=0$.

 Let now $\g =sl_2(\C) $ with generators $e$, $f$, $h$ and
relations $[h,e]= 2 e$, $[h, f] = -2 f$, $[e,f]= h$.
Let $\Lambda_0$, $\Lambda_1$  be the fundamental weights for $\hg$.

For $s \in {\Z}$, we define $H^{s} = - s \frac{h}{2}$. Then $(
H^{s} , h ) = -s$.
Define $$\Delta_s (z) = z^{ H^{s}(0)} \exp \left(\sum_{ n=1}
^{\infty} \frac{H^{s}(n)}{-n} (-z)^{-n}\right).$$
Applying  the results obtained in  \cite{Li5} on $L(k\Lambda_0)$--modules
we get the following proposition.

\begin{proposition} \label{novi-op} Let $s \in {\Z}$ and let $\Delta_s (z)$ be
defined as above.
For any weak $L(k\Lambda_0)$--module $(M,Y_M(\cdot,
 z))$,
$$ (\pi_s(M), Y_M ^{s}(\cdot,z)):=(M, Y_M (\Delta_s (z) \cdot, z))$$
is a weak $L(k \Lambda_0)$--module. $\pi_s(M)$ is an irreducible weak $L(k\Lambda_0)$--module if and only if
  $M$ is   an irreducible weak $L(k\Lambda_0)$--module.
\end{proposition}

By definition we have:
\bea &&  \Delta_s (z) e(-1){\1} = z^{-s} e(-1){\1}, \nonumber \\
 &&
\Delta_s (z) f(-1){\1} = z^{s} f(-1){\1}, \nonumber  \\
&& \Delta_s (z) h(-1){\1} = h(-1){\1} - s k z^{-1} {\1}. \nonumber
\eea

In other words, the corresponding automorphism $\pi_s$ of $U(\hg)$
satisfies the condition:
$$ \pi_s ( e(n) ) = e(n - s), \ \ \pi_s ( f (n) ) = f (n+s), \ \
\pi_s( h(n) ) = h(n) - s k \delta_{n,0}.$$

In the case $s=-1$ we get $$\pi_{-1}( L(k\Lambda_0)) = L(k
\Lambda_1).
$$
It is also important to notice the following important property:
$$ \pi_{s + t} (M) \cong \pi_s ( \pi_t (M) ), \ \ (s, t \in
{\Z}).$$
In particular, $$ M \cong \pi_0 (M) \cong \pi _s ( \pi_{-s} (M) ).
$$

\begin{remark}
Modules of type $\pi_s(M)$ have played an important role in the fusion rules analysis (cf. \cite{Gab}).
\end{remark}

\section{Vertex operator algebra $L(-\frac{4}{3} \Lambda_0)$}

First  we   study the vertex operator algebra $N(-\frac{4}{3}
\Lambda_0)$. Since level $-\frac{4}{3}$ is admissible, this vertex
operator algebra contains a unique maximal ideal $N^{1}
(-\frac{4}{3} \Lambda_0)$ which is generated by one singular
vector $v_{sing}$ (cf. \cite{AM}, \cite{DLM}, \cite{Gab}, \cite{KW}).
This vector can be written in the PBW basis of $N(-\frac{4}{3}\Lambda_0)$ as follows :
$$ v_{sing} = f_{sing}{\1}$$
where $f_{sing} \in U(\hg)$ is defined by
$$ f_{sing} = e(-1) \phi_{-\frac{4}{3}} + \frac{1}{3} e(-3)  - \frac{1}{2}
h(-1) e(-2) + \frac{1}{2} h(-2) e(-1),$$
and
$$ \phi_{-\frac{4}{3} } =  \frac{3}{4} ( e(-1) f(-1) + f(-1) e(-1) +
\frac{1}{2} h(-1) ^{2} ). $$
Note that $\omega_{Sug} =\phi_{-\frac{4}{3} } {\1}$ is the usual
Virasoro element in the vertex operator algebra $N(-\frac{4}{3}
\Lambda_0)$ obtained by  using the Sugawara construction (cf.
\cite{K2}, \cite{FZ}, \cite{MP}).

Set $$L(-\frac{4}{3} \Lambda_0)= \frac{N(-\frac{4}{3}
\Lambda_0)}{N ^{1}(-\frac{4}{3} \Lambda_0)}. $$
Then $L(-\frac{4}{3} \Lambda_0)$ is a simple vertex operator
algebra of rank $-6$.  The representation theory of this vertex
operator algebra is very interested.

 It is  important  to determine  which ${\hg}$--modules of level
 $-\frac{4}{3}$ have the structure of a module for the vertex
 operator algebra $L(-\frac{4}{3} \Lambda_0)$. This problem
 can be solved inside certain categories of $\hg$--modules (see \cite{AM}, \cite{DLM}).

 Applying the general result obtained in   Theorem 3.5.3 of \cite{AM}
 to  the case when  $k=-\frac{4}{3}$
 we get the following result.

\begin{theorem} \cite{AM} \label{klas-0}
The set
$$ \{ L(-\frac{4}{3} \Lambda_0), \ L(-\frac{2}{3} \Lambda_0 - \frac{2}{3}
\Lambda_1), \ L(-\frac{4}{3} \Lambda_1) \}
$$
provides all irreducible $L(-\frac{4}{3} \Lambda_0)$--modules from
the category ${\mathcal O}$. Every $ L(-\frac{4}{3}
\Lambda_0)$--module from the category ${\mathcal O}$ is completely
reducible.
\end{theorem}

Applying the automorphism $\pi_s$, $s \in {\Z}$, on
${\ver}$--modules from the category $\mathcal{O}$ one can get a
family of weak ${\ver}$--modules which in general don't belong to
the category $\mathcal{O}$. These modules can be
identified by using the following lemma.

\begin{lemma} \label{sing-za}
\item[(1)] Assume that $M$ is a weak ${\ver}$--module which is
generated by the vector $v_{s}$ $(s \in {\Z})$ such that:
\bea && e(n- s) v_{s} = f (n+s) v_{s}=0 \quad (n \ge 0),
\label{rel-0-1} \\
&& h(n) v_{s} = \delta_{n,0} ( - \frac{4}{3} s ) v_{s} \quad (n
\ge 0).\label{rel-0-2} \eea
Then
$$ M \cong \pi_{-s} ( L( -\frac{4}{3}
\Lambda_0)). $$
\item[(2)] Assume that $M$ is a weak ${\ver}$--module which is
generated by the vector $v'_{s}$ $(s \in {\Z})$ such that:
\bea && e(n- s) v'_{s} = f (n+1+s) v'_{s}=0 \quad (n \ge 0),
\label{rel-1-1} \\
&& h(n) v'_{s} = \delta_{n,0} ( -\frac{2}{3}- \frac{4}{3} s )
v'_{s} \quad (n \ge 0).\label{rel-1-2} \eea
Then
$$ M \cong \pi_{-s} ( L( -\frac{2}{3}
\Lambda_0 - \frac{2}{3} \Lambda_1)). $$
\end{lemma}
{\em Proof.}  (1) We consider the weak $L(-\frac{4}{3}
\Lambda_0)$--module $\pi_{s} (M)$. By construction, we have that
$\pi_{s}(M)$ is a highest weight $\hg$--module with the highest
weight $-\frac{4}{3}  \Lambda_0 $. Since every ${\ver}$--module
from the category $\mathcal{O}$ is completely reducible, we have
that
$\pi_{s}(M) \cong L(-\frac{4}{3} \Lambda_0),$
which implies that
$$ M \cong \pi_{-s} ( \pi_{s} (M) ) \cong \pi_{-s} (L(-\frac{4}{3}   \Lambda_0)  ),$$
and the statement (1)  holds. The proof of (2) is similar.  \qed

\vskip 5mm

\section{Lattice construction of $L(-\frac{4}{3} \Lambda_0)$ and its modules}
\label{lattice-construction}

 In this section we shall present a lattice construction of the
vertex operator algebra $L(-\frac{4}{3} \Lambda_0)$.
This construction will use certain ideas and  results
obtained in \cite{A}.
 Our results in  Section \ref{podalgebre} will show that this approach can be applied because
  the vertex operator algebra $L(-\frac{4}{3}
\Lambda_0)$ contains the $\mathcal{W}(2,5)$--algebra investigated in \cite{A}.
\vskip 5mm
Define first the following lattice:
$$ \widetilde{L}= {\Z}\a + {\Z} \b, \qquad \la \a , \a \ra = - \la \b , \b \ra = \frac{1}{6}, \ \la \a , \b \ra = 0.$$

Let
 ${\h}={\C}\otimes_{\Z} \widetilde{L}$.
  Extend the form $ \la \cdot,
\cdot \ra $ on $\widetilde{L}$ to ${\h}$.
 Let $\hat{{\h}}={\h} \,\otimes{\C}[t,t^{-1}] \oplus {\C}c$ be the affinization of
${\h}.$
Set $ \hat{{\h}}^{+}={\h} \otimes t{\C}[t]
;\;\;\hat{{\h}}^{-}={\h} \otimes t^{-1}{\C}[t^{-1}]. $ Then
$\hat{{\h}}^{+}$ and $\hat{{\h}}^{-}$ are abelian subalgebras of
$\hat{{\h}}$. Let $U(\hat{{\h}}^{-})=S(\hat{{\h}}^{-})$ be the
universal enveloping algebra of $\hat{{\h}}^{-}$. Let ${\l} \in
{\h}$. Consider the induced $\hat{{\h}}$-module
\begin{eqnarray*}
M(1,{\l})=U(\hat{{\h}})\otimes _{U({\h} \otimes {\C}[t]\oplus
{\C}c)}{\C}_{\l}\simeq
S(\hat{{\h}}^{-})\;\;\mbox{(linearly)},\end{eqnarray*} where ${\h}
\otimes t{\C}[t]$ acts trivially on ${\C}$,
${\h}$ acting as $\la \alpha , {\l} \ra$ for $\alpha \in {\h}$
and $c$ acts on ${\C}$ as multiplication by 1. We shall write
$M(1)$ for $M(1,0)$.
 For $\alpha \in {\h}$ and $n \in {\Z}$ write $\alpha (n) =
 \alpha\otimes t^{n} $. Set
$ \alpha(z)=\sum _{n\in {\Z}}\alpha(n)z^{-n-1}. $
Then $M(1)$ is a vertex   algebra which is generated by the fields
$\alpha(z)$, $\alpha \in {\h}$, and $M(1,{\l})$, for $\l \in
{\h}$, are irreducible modules for $M(1)$.

As in  \cite{DL}, \cite{GL} (see also \cite{FLM}, \cite{K2}), we
have the generalized vertex algebra
$$ V_{\widetilde{L}}  = M(1) \otimes {\C}[\widetilde{L}], $$
where ${\C}[\widetilde{L}]$ is a group algebra of $\widetilde{L}$
with generators $e ^{\a}$ and $ e^{\b}$.
For $v \in  V_{ \widetilde{L} }$ let $Y(v,z) = \sum_{ s \in
\frac{1}{6} {\Z} } v_s z^{-s -1}$ be the corresponding vertex
operator (for precise formulae see \cite{DL}).

Recall some useful facts which hold in the generalized vertex
algebra $V_{\widetilde{L}}$.

Assume that $\alpha^{1}, \alpha^{2} \in \widetilde{L}$, $\la
\alpha^{1} , \alpha^{2} \ra = r \in {\Z}$. Then
\bea
&& Y(e^{\alpha ^{1}}, z) e^{\alpha^{2}} = \sum_{n \in {\Z} }
e^{\alpha^{1}}_{n+r} e^{\alpha^{2}} z^{-n-r-1}, \label{rel-gen-0} \\
&& e^{\alpha^{1}}_{ - r + n} e^{\alpha^{2}} = 0 \ \ \ \mbox{for every} \ n
\in {\Zp}. \label{rel-gen-2} \eea

Clearly, $M(1)$ is a vertex subalgebra of $V_{ \widetilde{L} }$.

Let $M_{\a}(1)$ (resp. $M_{\b}(1)$)  be a vertex subalgebra of
$M(1)$ generated by the field $\a  (z)$ (resp. $\b(z)$).
$M_{\a}(1)$ (resp. $M_{\b}(1)$)  is a subalgebra of the
generalized vertex algebra $V_{\Z \a}$ (resp. $V_{\Z \b}$).

Define the following Virasoro element in $M_{\a}(1) \subset V_{\Z
\a} \subset V_{\widetilde{L}}$:
\bea \label{exp-vir} && \omega ^{\a}  = 3 \a  (-1) ^{2} - 2 \a (-2) . \eea
Then $\omega^{\a}$ spans a subalgebra of  $M_{\a}(1)$ isomorphic
to the Virasoro vertex operator algebra $L^{Vir} (-7,0)$ with
central charge $-7$.

As in \cite{A}, we define the following screening operators:
\bea \label{exp-scr} Q = e^{-6 \a}_0, \quad \widetilde{Q} = e^{2 \a} _0, \eea
and the following subalgebras of $M(1)$:
$$\overline{M(1)} = \mbox{Ker} _{M(1)} \widetilde{Q}, \quad
\overline{M_{\a}(1)} = \mbox{Ker} _{M_{\a} (1)} \widetilde{Q}.$$
Then
\bea \label{exp-h} H= e^{-6 \a} _0 e^{6 \a } \in \overline{M_{\a}(1)} \subset \overline{M(1)} \eea
is a primary vector of conformal weight $5$.
\begin{remark} \label{napomena-1}
The generalized vertex algebra $V_{\widetilde{L}}$ is larger than
the algebras considered in \cite{A} in the   case when $p=3$. More
precisely,  if we set ${\ap}= - 6 {\a}$ and ${\bp}= - 2 {\a}$,
then the subalgebras $V_{\Z {\ap} }$ and $V_{\Z {\bp}}$ of
$V_{\widetilde{L}}$ are isomorphic to the (generalized) vertex
operator algebras considered in \cite{A}. Since
$$ Q=e^{\ap}_0, \quad \widetilde{Q}= e^{-\bp}_0, \quad \omega ^{\a}= \frac{3}{4} {\bp} (-1) ^{2} + {\bp}(-2),
\quad H= e^{\ap}_0 e^{-\ap}, $$
we have that
 expressions for  vectors and operators defined by relations (\ref{exp-vir})-(\ref{exp-h}) coincide with
 the definition of the same  objects in Section 2 of \cite{A} when $p=3$.
\end{remark}

The results from \cite{A} applied to the case $p=3$ give the following
result:
\begin{theorem} \cite{A} \label{res-a}
\item[(1)] $[Q, \widetilde{Q}] = 0$,
\item[(2)] $\mbox{Ker} _{ M_{\a} (1)} Q \cong L^{Vir} (-7,0)$,
\item[(3)]
The vertex operator algebra $ \overline{M_{\a}(1)}$ is generated
by $\omega^{\a}$ and $H= Q e^{6 \a}$.
\end{theorem}

This theorem has the following consequence:

\begin{corollary} \label{kor-a}
\item[(1)] The vertex operator algebra $\mbox{Ker} _{ M(1) } Q
\cong L^{Vir}(-7,0) \otimes M_{\b}(1)$ is generated by $\b (-1)$
and $\omega^{\a}$.
\item[(2)]The vertex operator algebra $\overline{M(1)} \cong
\overline{M_{\a}(1)} \otimes M_{\b}(1) $   is generated by $\b
(-1)$, $\omega^{\a}$ and $H$.
\end{corollary}
Define
$$ \omega = \omega^{\a} - 3 \b(-1) ^{2} = 3 \a(-1)^{2} - 2 \a (-2) -
3 \b(-1)^{2}. $$
Then $\omega$ is a Virasoro element which generates the subalgebra
$L^{Vir}(-6,0)$. Set
$$L(z)= Y(\omega,z) = \sum_{ n \in {\Z} } L(n) z^{-n-2}. $$

Define the following vectors in $ V_{\widetilde{L}}$:
\bea
& e &= e^{3 (\a - \b)}, \label{def-e} \\
& h & = 4 \b (-1) , \label{def-h} \\
& f & = -\frac{2}{9} Q e^{3 (\a + \b)} = - ( 4 \a(-1) ^{2} -
\frac{2}{3}\a(-2) ) e^{-3(\a - \b)}. \label{def-f}
\eea
Then $e, h, f$ are primary vectors of conformal weight $1$ for the
Virasoro algebra, i.e.,
$$
L(n) e = \delta_{n,0} e, \quad L(n) h = \delta_{n,0} h, \quad L(n)
f = \delta_{n,0} f \quad (n \ge 0). $$
We are interested in the subalgebra generated by $e$, $f$ and $h$.
For this purpose we define the following lattice:
$$  {L}= {\Z}(\a + \b) + {\Z}(\a - \b).$$
Clearly, $V_{ {L}}$ is a subalgebra of $V_{\widetilde{L}}$ which
contains vectors $e$, $f$ and $h$. It is important to notice that
\bea \label{vazna}&& \la 3(\a \pm \b),  {L} \ra  \subset {\Z},
    \eea
i.e., $3(\a + \b)$ and $3(\a - \b)$ are elements of the dual
lattice of $L$ (see also Remark \ref{dualna} below).

The following lemma is a  consequence of (\ref{vazna}) and the
Jacobi identity in  generalized vertex algebras (see \cite{DL},
\cite{GL}, \cite{S}).

\begin{lemma} \label{lema-komutator}
Assume that $a, b \in \mbox{span}_{\C} \{e, f, h \}$ and $ v \in
V_{ {L}}$. Then
$ Y( a, z) v = \sum_{n \in {\Z} } a_n v  z^{-n-1} $,
and for every $m, n \in {\Z}$ the following commutator formula
holds:
\bea \label{komutiranje} &&
[a_n, b_m] v = \sum_{i = 0} ^{\infty} {n \choose i} (a_i b)_{n +m
- i} v .  \eea
\end{lemma}

For $a \in \mbox{span}_{\C} \{e, f, h \}$ set $a(n)=a_n $ and
$a(z)= \sum_{n \in \Z} a(n) z^{-n-1}$. Lemma \ref{lema-komutator}
shows that the field $a(z)$ is in fact the restriction of the
vertex operator $Y(a,z)$ on $V_{{L}}$.

\begin{theorem} \label{teorem-podalgebra}
  The vectors $e$, $f$ and $h$ span a subalgebra of the
generalized vertex algebra $V_{{L}}$ isomorphic to the vertex
operator algebra ${\ver}$. Moreover, $V_{{L}}$ is a weak module
for the vertex operator algebra ${\ver}$.
\end{theorem}
{\em Proof.}
  First we notice that for
every $n \in {\Z}$,  $n \ge 0$, the following relations in the
generalized vertex algebra $V_{{L}}$ hold:
\bea
&e(n) e& = 0  , \nonumber \\
&f(n) f& = \frac{2}{81} Q^{2} ( e^{3(\a + \b)} _n e^{3(\a + \b)} ) = 0 ,   \nonumber \\
&h(n) h & = 0 \quad   (n \ne 1),
\nonumber \\
& h(1) h &= -\frac{8}{3} {\1}, \nonumber \\
& h(n) e &= 2 \delta_{n, 0} e, \nonumber \\
& h(n) f &= - 2 \delta_{n, 0} f, \nonumber \\
& e(n) f &= 0 \quad \mbox{for every} \ n \ge 2, \nonumber \\
& e(1) f &=  - \frac{4}{3} {\1}, \nonumber \\
& e(0) f &= h . \nonumber
 \eea
By using these relations and  commutator formula
(\ref{komutiranje}) we conclude that the components of the fields
$e(z), f(z), h(z)$  satisfy the commutation relations for the
affine Lie algebra $\hg$ of level $  -\frac{4}{3}$. So $V_{{L} }$
is a $U(\hg)$--module of level $-\frac{4}{3}$.

 By construction, we
have that $V_{{L} }$ is a weak module for the vertex operator
algebra $N(-\frac{4}{3}\Lambda_0)$.
Let $S$ be the subalgebra of $V_L$  generated by  vectors $e$, $f$
and $h$, i.e.,
$$  S=\mbox{span}_{\C} \{ u^{1}_{n_1} \cdots u^{r} _{n_r}
{\1}  \vert \   u^{1}, \dots, u^{r} \in \{e,f,h\}, \ n_1, \dots,
n_r \in {\Z}, r \in {\Zp} \}. $$
Then $S$ is isomorphic to a certain quotient of the vertex
operator algebra $N(-\frac{4}{3} \Lambda_0)$.
 The Virasoro element $\omega$ belongs to the subalgebra $S$ because
$$\phi_{-\frac{4}{3}}.{\1}=\omega = 3 \a (-1) ^{2} - 2 \a (-2) - 3
\b(-1)^{2}. $$
 As a
$U(\hg)$--module, $S$  a cyclic  module generated by the vacuum
vector ${\1}$.
 Next we consider element $f_{sing}. {\1} \in S$. The definition of the
 action of $\hg$ on $V_L$ gives that
 \bea
  f_{sing} . {\1} &=& e(-1) \omega +\frac{1}{3} e(-3){\1} -
  \frac{1}{2}h(-1) e(-2){\1} + \frac{1}{2} h(-2) e(-1){\1}
  \nonumber \\
  & = & e^{3 (\a - \b)} _{-1} \omega + \frac{1}{3} e^{3 (\a - \b)} _{-3} {\1} - 2 \b(-1) e ^{3(\a - \b)}_{-2} {\1}
  + 2 \b(-2) e^{3 ( \a -  \b)} \nonumber \\
& = & \frac{4}{3} e^{3(\a - \b)} _{-3} {\1} - 9 ( \a(-1) - \b(-1)
)^{2} e^{3 (\a - \b)}\nonumber \\ & & - 6 \b(-1) ( \a(-1) -
\b(-1)) e^{3(\a - \b)}  + 2 \b(-2) e^{3(\a - \b)}  \nonumber \\ &
& + ( 3 \a(-1)^{2} - 2 \a(-2) - 3
\b(-1)^{2} ) e^{3(\a - \b)} \nonumber \\
&=&  \frac{2}{3}\left( (3 \a(-1) - 3 \b(-1) )^{2} + 3 \a(-2) - 3
\b(-2) \right) e^{3(\a - \b)} \nonumber \\
& & + ( - 9 ( \a(-1) - \b(-1) )^{2}  - 6 \b(-1) ( \a(-1) - \b(-1))
\nonumber \\
& &
 + 2 \b(-2)  + 3 \a(-1)^{2} - 2
\a(-2) - 3
\b(-1)^{2} ) e^{3(\a - \b)} \nonumber \\
& = & 0. \nonumber
 \eea
Since $f_{sing} . {\1} = 0$, we conclude that $S$ is a certain
quotient of $$\frac{N(-\frac{4}{3}\Lambda_0)}{  U(\hg) v_{sing} }
=L(-\frac{4}{3}\Lambda_0).$$
The simplicity of ${\ver}$ implies that $S=U(\hg).{\1} \cong
 {\ver}$.
  Therefore, $V_{{L}}$ is a weak ${\ver}$--module.
 \qed

\begin{remark} \label{dualna}
The vertex operator algebra  ${\ver}$ can be embedded into certain
vertex algebras. Let us describe these algebras. Define
$$ L_0=  {\Z} (3{\a} -3 {\b}) + {\Z}  (3{\a} + 3{\b}).$$
Then $L_0$ is the dual lattice of $L$ in $\h$.
  $L_0$ is an even lattice and
therefore the subalgebra $V_{L_0} $ of $V_{L}$ has the structure
of a vertex algebra (cf. \cite{DL}). Moreover, $V_L$ is a
$V_{L_0}$--module. Since $e,f, h \in V_{L_0} $, we have that the
vertex operator algebra $L(-\frac{4}{3} \Lambda_0)$ is a
subalgebra of $V_{L_0} $. It is also interesting to notice that
$ e, f, h \in \mbox{Ker} _ {V_{L_0}} \widetilde{Q}$.
So ${\ver}$ is also a subalgebra of the vertex algebra $\mbox{Ker}
_ {V_{L_0} } \widetilde{Q}$.

On the other hand, we also want to construct explicitly other
${\ver}$--modules from the category ${\mathcal O}$. It turns out
that for this purpose the larger algebra $V_{ {L}}$ is more
suitable.
\end{remark}

Now we want to identify certain ${\ver}$--submodules of $V_{L}$.
We have the following result.

\begin{theorem} \label{spectral} For every $s \in {\Z}$ we have:
\item[(1)] $$U(\hg) . e^{ 2 s \b} \cong \pi _{-s} (L(-\frac{4}{3}
\Lambda_0)).$$

\item[(2)]$$U(\hg) . e^{-\a + (2 s +1 )\b} \cong \pi _{-s}
(L(-\frac{2}{3} \Lambda_0 - \frac{2}{3} \Lambda_1)). $$
\end{theorem}
{\em Proof.}
Define
$v_s = e^{2 s \b}$ and $v'_s =  e^{-\a + (2 s +1 )\b}$.
Since $V_L$ is a weak ${\ver}$--module, we have that for every $s
\in {\Z}$ the submodules
$U(\hg). v_s$ and $U(\hg) . v'_s $
have the structure of a weak module for the vertex operator
algebra ${\ver}$.

Assume that $n \ge 0$. Since $ \la 4 \b , 2 s \b \ra =
-\frac{4}{3} s$, we have
$$h(n) v_s = -\frac{4}{3} s \delta_{n,0} v_s . $$
Since $\la 3 \a - 3\b , 2 s \b \ra = s$, relation (\ref{rel-gen-2})
implies that
$$ e( n - s) v_s = 0. $$
Since  $\la 3 \a  + 3\b , 2 s \b \ra = -s$ and $Q v_s = 0$, we
have that
\bea
f(n+s) v_s &=& -\frac{2}{9} (Q e^{3\a + 3 \b})_{n+s}   v_s
\nonumber \\
&=&-\frac{2}{9} Q e^{3\a + 3 \b}_{n+s}  v_s + \frac{2}{9} e^{3\a
+ 3 \b}_{n+s}  Q v_s \nonumber \\
& =& 0 . \nonumber \eea
In this way we have verified that vector $v_s$ satisfies
conditions (\ref{rel-0-1}) and (\ref{rel-0-2}) and Lemma
\ref{sing-za} implies that
$$ U(\hg) . v_s \cong \pi_{-s} (L(-\frac{4}{3} \Lambda_0)).$$
This proves assertion (1).

%
%
 Since $ \la 4 \b , -\a + (2 s + 1) \b \ra =
-\frac{2}{3}-\frac{4}{3} s$, we have
$$h(n) v'_s =(-\frac{2}{3} -\frac{4}{3} s ) \delta_{n,0} v'_s . $$
Since $\la 3 \a - 3\b , - \a + ( 2 s + 1) \b \ra = s$, relation
(\ref{rel-gen-2}) implies that
$$ e( n - s) v'_s = 0. $$
Since  $\la 3 \a  + 3\b , -\a +(2 s +1)\b \ra = -s-1$ and $Q v'_s
= 0$, we have that
\bea
f(n+1+s) v'_s &=& -\frac{2}{9} (Q e^{3\a + 3 \b})_{n+1+s}   v'_s
\nonumber \\
&=&-\frac{2}{9} Q e^{3\a + 3 \b}_{n+1+s}  v'_s + \frac{2}{9}
e^{3\a
+ 3 \b}_{n+1+s}  Q v'_s \nonumber \\
& =& 0 . \nonumber \eea
So vector $v'_s$ satisfies conditions (\ref{rel-1-1}) and
(\ref{rel-1-2}) and by Lemma \ref{sing-za} we have that
$$ U(\hg) . v'_s \cong \pi_{-s} (L(-\frac{2}{3} \Lambda_0 - \frac{2}{3}\Lambda_1)),$$
and assertion (2) holds. \qed

 This theorem has the following consequence.

\begin{corollary} The vectors
$e^{-\a + \b}$, $e^{2\b}$ are singular vectors for the action of
$\hg$. Moreover, we have
$$ U(\hg) . e^{-\a + \b} \cong  L(-\frac{2}{3} \Lambda_0 -
\frac{2}{3}\Lambda_1), \quad U(\hg) . e^{2 \b} \cong L(-
\frac{4}{3}\Lambda_1). $$
\end{corollary}

The generalized vertex algebra $V_L$ contains some other interesting ${\ver}$--submodules. Let us give one example.

\begin{example}
Define
$ \widetilde{E} = U(\hg) . e ^{-2 \a}$.
Since
$ {\g} \otimes t {\C} [t] . e^{-2 \a} =0$,
we have that the $\hg$--module $\widetilde{E}$ is
${\Zp}$--graded.
Moreover, the top level
 $ E= U(\g) . e^{-2 \a}$
is an irreducible $U(\g)$--module which is neither highest nor lowest weight with respect to ${\g}$.
Therefore, the $\hg$--module $\widetilde{E}$ is not in the category
$\mathcal{O}$, but it is a ${\Zp}$--graded $L(-\frac{4}{3}
\Lambda_0)$--module  on which $h(0)$ and $L(0)$ act semisimply.
By using Zhu's algebra theory (cf. \cite{Z}),  modules of this type were
investigated in \cite{AM}.
\end{example}

\section{Vertex subalgebras of $L(-\frac{4}{3}
\Lambda_0)$}
\label{podalgebre}

In Section \ref{lattice-construction} we proved that the vertex
algebra $M(1)$ contains the subalgebra $\overline{M(1)}$ which is
a vertex operator algebra of rank $-6$. This vertex operator
algebra is obtained as kernel of the screening operator
$\widetilde{Q}$. In this section we will show that
$\overline{M(1)}$ actually lies in the vertex operator algebra
${\ver} \subset V_L$. This result will imply that the coset vertex
operator algebra
$$C({{\ver} }, M_{\b}(1) )\cong
\overline{M_{\a}(1)}.$$
Since $\overline{M_{\a}(1)}$ is a simple vertex operator algebra
 generated by Virasoro element $\omega^{\a}$ and   primary
vector $H$ of conformal weight $5$, we can say that this coset
vertex operator algebra (or the coset
  $\frac{\widehat{su(2)}_{-\frac{4}{3}} }{ \widehat{u(1)} } $ in the  terminology of \cite{BE})
is   the $\mathcal{W}(2,5)$ algebra with central charge $-7$.

\vskip 10mm

The operator $h(0)$ acts semisimply on ${\ver}$ and
$$ {\ver} = \bigoplus _{ i \in {\Z} } M ^{(i)} , \
\mbox{where} \ \   M^{(i)} = \{ v \in {\ver} \ \vert \ h(0) v = 2
i  v \}.$$

\begin{proposition} \label{direct-sum}
The vertex operator subalgebra $M^{(0)}=\mbox{Ker}_{\ver} h(0)$ is
simple. As a $M^{(0)}$--module, ${\ver}$ is isomorphic to the
direct sum $ {\ver} = \bigoplus _{ i \in {\Z} } M ^{(i)}$ , and
each $M^{(i)}$ is a simple $M^{(0)}$--module.
\end{proposition}
{\em Proof.} Clearly, $M^{(0)}=\mbox{Ker}_{\ver} h(0)$ is a vertex
operator subalgebra of ${\ver}$. For  $u \in M^{(i)}$ and $n \in
{\Z}$ we have that $u_n M^{(j)} \subset M^{(i+j)}$.

Let $0 \ne v \in M^{(j)}$.  Since ${\ver}$ is a simple vertex
operator algebra, by Corollary 4.2 of \cite{DM-galois}  we have
that
\bea \label{generiranje} && {\ver} = \mbox{span}_{\C} \{ u_n v \
\vert \ u \in {\ver}, \ n \in {\Z} \} . \eea
 This relation implies that
\bea \label{generiranje-2} && M^{(j)} = \mbox{span}_{\C} \{ u_n v
\ \vert \ u \in M^{(0)}, \ n \in {\Z} \} . \eea
 So $M^{(j)}$
is a simple $M^{(0)}$--module. In particular, $M^{(0)}$ is a
simple vertex operator algebra. \qed

\vskip 5mm

We shall now identify the subalgebra $M^{(0)}$.

\begin{theorem} \label{izom-w}
The vertex operator algebra   $M^{(0)}= \mbox{Ker}_{ {\ver}} h(0)$
is isomorphic to $\overline{M(1)} \cong \overline{M_{\a}(1)}
\otimes M_{\b}(1)$. In particular, $\overline{M(1)}$ and
$\overline{M_{\a}(1)}$ are simple vertex operator algebras.
\end{theorem}
{\em Proof.}
 By using the lattice construction of
${\ver}$ from Section \ref{lattice-construction}, one can easily
see that $M^{(0)}$ is actually a subalgebra of the vertex operator
algebra $M(1)$. Since ${\ver} \subset \mbox{Ker}_{V_L}
\widetilde{Q}$, we have that
\bea \label{inkluzija-1} && M^{(0)} \subseteq \mbox{Ker}_{ M(1)} \widetilde{Q} = \overline{M(1)} .
\eea
We shall now  prove the other inclusion. By   Corollary \ref{kor-a} we have that
the vertex operator algebra $\overline{M(1)}$ is generated by
vectors $\b(-1) $, $\omega^{\a}$ and $H$. It suffices to prove
that these generators belong to ${\ver}$. Clearly,
\bea \label{beta} && {\b} (-1) = \frac{1}{4} h \in {\ver}. \eea
The fact that  $\omega = \omega^{\a} - 3 {\b}(-1) ^{2}$ is the
Virasoro element in ${\ver}$ implies that
\bea && \label{omega-a} \omega^{\a} = \omega + \frac{3}{16} h(-1) ^{2} {\1}  \in {\ver}. \eea
Since $\mbox{Ker}_{M(1)} Q $ is generated by $\b(-1)$ and
$\omega^{\a}$, relations (\ref{beta}) and (\ref{omega-a}) imply
that
\bea \label{ul-q} && \mbox{Ker}_{M(1)} Q \subset M^{(0)}. \eea
Now we shall prove that $H \in M^{(0)}$. Define vector
$$ v= -9 e(-4) f = 2 e^{3\a - 3 \b} _{-4} Q e^{3\a + 3\b} \in M^{(0)}.$$
Since $ Q^{2} e^{3\a - 3 \b} = Q^{2} e^{3\a + 3 \b} = 0$, we have
$$ Q v =  2 ( Q e^{3\a - 3 \b}  )_{-4} (Q e^{3\a + 3 \b}) = Q^{2} (
e^{3\a - 3\b} _{-4} e^{3\a + 3 \b}  ) = Q^{2} e^{6 \a} . $$
This implies that
$ Q( v - Q e^{6\a}) = 0$.
Since
$$v - Q e^{6\a} \in \mbox{Ker}_{M(1)} Q \subset M^{(0)}$$
we conclude that
\bea \label{h} &&H = Q e^{6\a} \in M^{(0)}. \eea
Now relation (\ref{inkluzija-1})-(\ref{h}) imply that
$M^{(0)}=\overline{M(1)}$.

The simplicity of $\overline{M_{\a}(1)}$ and $\overline{M(1)}$
follows from the fact that $M^{(0)}$ is a simple vertex operator
algebra (cf. Proposition \ref{direct-sum}).\qed
\vskip 5mm

By the lattice construction, we have that the subalgebra of
${\ver}$ generated by the field $h(z)$ is isomorphic to
$M_{\b}(1)$. Theorem \ref{izom-w} enables us to describe the
following vertex operator algebra.

\begin{corollary}
The coset vertex operator algebra
 $$C({{\ver} }, M_{\b}(1) ) =\{ v \in {\ver} \ \vert \ h(n) v = 0, \ n \ge 0 \}$$
is isomorphic to the vertex operator algebra
$\overline{M_{\a}(1)}$ (=$\mathcal{W}(2,5)$ algebra with central charge $-7$).
\end{corollary}

\begin{remark}
It was proved in  \cite{KR} that the subalgebra
$\mbox{Ker}_{L(-\frac{1}{2}\Lambda_0)} h(0)$ is isomorphic to the
vertex algebra $\mathcal{W}_{1+\infty}$ with central charge $c=-1$. In
\cite{W2}, Wang proved  that the vertex algebra $\mathcal{W}_{1+\infty}$
with central charge $-1$ is isomorphic to the tensor product of
the $\mathcal{W}(2,3)$ algebra with central charge $-2$ and the free boson
vertex algebra. Our result shows that the vertex operator algebra
${\ver}$ contains a subalgebra of similar type. The difference is
that the $\mathcal{W}(2,3)$ algebra with central charge $-2$ is replaced
by  the $\mathcal{W}(2,5)$ algebra with central charge $-7$.
\end{remark}

\end{document}